\documentclass[leqno,12pt]{article}  
\usepackage{amsmath}
\usepackage{amssymb}

\usepackage[german,english]{babel}
\usepackage{amsthm}
\usepackage{xypic}  

\author{\textsc{Elmar Grosse-Kl\"onne}}
\date{}

\theoremstyle{plain} 
\newtheorem{satz}{Theorem}[section]  
\newtheorem{kor}[satz]{Corollary}  
\newtheorem{lem}[satz]{Lemma}  
\newtheorem{pro}[satz]{Proposition}  
 
\newcommand{\spm}{\mbox{\rm Sp}}  
\newcommand{\spf}{\mbox{\rm Spf}}  

\newcommand{\codim}{\mbox{\rm codim}}  
\newcommand{\bi}{\mbox{\rm Im}}  
\newcommand{\ke}{\mbox{\rm Ker}}  

\newcommand{\q}{\mbox{\rm Frac}}  
\newcommand{\kara}{\mbox{\rm char}}  

\newcommand{\kodi}{\mbox{\rm coh.dim}}

\theoremstyle{remark}

\theoremstyle{definition}

\begin{document}

%

\begin{center}{\bf De Rham Cohomology Of Rigid Spaces}\\Elmar Grosse-Kl\"onne
\end{center}
\begin{abstract} We define de Rham cohomology groups for rigid spaces over non-archimedean fields of characteristic zero, based on the notion of dagger space introduced in \cite{en1dag}. We establish some functorial properties and a finiteness result, and discuss the relation to the rigid cohomology as defined by P. Berthelot \cite{berco}. 
\end{abstract}

\begin{center}{\bf Introduction}\end{center}

Let $k$ be a field, complete with respect to a non-archimedean valuation. V. Berkovich, R. Huber and others have worked out satisfying foundations for the {\it \'{e}tale} cohomology theory of $k$-rigid spaces. In this paper we assume $\kara(k)=0$ and propose a definition of {\it de Rham} cohomology groups $H^*_{dR}(X)$ for $k$-rigid analytic spaces $X$. In \cite{ss} it is shown that for smooth $X$ the naive definition $H^*(X,\Omega^{\bullet}_X)$ meets certain minimal axiomatic requirements. However, if $X$ presents "boundaries" there are also serious pathologies entailed by this definition. For example, if $X=\spm(T_1)$, the closed unit disk, then $H^1(X,\Omega_X^{\bullet})$ is an infinite dimensional $k$-vector space; formal integration preserves the radius of convergence, but does not preserve convergence on the boundary. 
The classical idea from the paper \cite{mw} of Monsky and Washnitzer to remedy this is to use only overconvergent power series for the definition of $H^*_{dR}(X)$. In \cite{en1dag} we introduced a category of 'rigid spaces with overconvergent structure sheaf' which we called $k$-dagger spaces, and in \cite{en2dag} we proved that $H^*(Y,\Omega^{\bullet}_Y)$ for many smooth $k$-dagger spaces $Y$ is indeed finite dimensional. Moreover, in \cite{en1dag} we associated to a $k$-dagger space functorially a $k$-rigid space with the same underlying $G$-topological space. 
Therefore, ideally we should define $H^*_{dR}(X)$ for a smooth $k$-rigid analytic space $X$ as $H^*(Y,\Omega^{\bullet}_Y)$, where $Y$ is a smooth $k$-dagger space whose associated rigid space is $X$. (For affinoid $k$-rigid spaces $X$, for which such a $Y$ exists, this is the approach of M. van der Put \cite{vdp}). Our results here imply that this definition is indeed independent, up to canonical isomorphism, of the choice of $Y$ 
--- if such a $Y$ exists. But if $X$ is not smooth, or if $X$ is not associated with a dagger space $Y$, another definition is needed. We propose: Let $W$ be a smooth $k$-dagger space and let $\phi:X\to W'$ be a closed embedding into its associated $k$-rigid space $W'$. Let $w:W'\to W$ be the natural morphism of ringed spaces (actually an isomorphism of underlying Grothendieck topologies). Then we set $$H^*_{dR}(X)=H^*(X,(w\circ\phi)^{-1}\Omega^{\bullet}_{W/k}).$$The content of the paper is organized as follows. In section 0 we recall some facts on dagger spaces, and in section 1 we show that $H^*_{dR}(X)$ is well-defined. In section 2 we establish some functorial properties, similar to \cite{hadr}. In section 3 we specialize to the case where $k=\q(R)$ with a complete discrete valuation ring $R$ of mixed characteristic, with residue field $\bar{k}$. From \cite{en2dag} we derive the finiteness of $H^*_{dR}(X)$ for a big class of $k$-rigid spaces $X$. We explain how the rigid cohomology, defined by Berthelot in \cite{berco}, of a finite type $\bar{k}$-scheme $Y$ can be expressed (or: be redefined) as the de Rham cohomology of the tube of $Y$ in a smooth formal $R$-scheme with an embedding of $Y$. A Gysin sequence for the de Rham cohomology of tubes in a semi-stable formal $R$-scheme is presented, generalizing the Gysin sequences for rigid cohomology from \cite{berfi}, \cite{mebend}. In section 4 we briefly discuss for a given admissible proper formal $R$-scheme the relation between the rigid cohomology of its special fibre and the de Rham cohomology of its generic fibre.

\addtocounter{section}{-1}
\section{Dagger spaces}

Let $k$ be a field complete with respect to a non-trivial non-archimedean valuation $|.|$, and of characteristic zero. We denote by $k_a$ its algebraic closure with value group $\Gamma^*=|k_a^*|=|k^*|\otimes\mathbb{Q}$.\\ 
We gather some facts from \cite{en1dag}. For $\rho\in\Gamma^*$ the $k$-affinoid algebra $T_n(\rho)$ consists of all series $\sum a_{\nu}X^{\nu}\in
k[[X_1,\ldots,X_n]]$ such that $|a_{\nu}|\rho^{|\nu|}$ tends to zero if $|\nu|\to\infty$. The algebra $W_n$ is defined to be $W_n=\cup_{\stackrel{\rho>1}{\rho\in\Gamma^*}}T_n(\rho)$. A $k$-dagger algebra $A$ is a quotient of some $W_n$; a surjection $W_n\to A$ endows it with a norm which is the quotient of the Gauss norm on $W_n$. All $k$-algebra morphisms between $k$-dagger algebras are continuous with respect to these norms, and the completion of a $k$-dagger algebra $A$ is a $k$-affinoid algebra $A'$ in the sense of \cite{bgr}. There is a tensor product $\otimes^{\dagger}_k$ in the category of $k$-dagger algebras.
As for $k$-affinoid algebras, one has for the set $\spm(A)$ of maximal ideals of $A$ the notions of rational and affinoid subdomains, and for these the analogue of Tate's acyclicity theorem (\cite{bgr},8.2.1) holds. The natural map $\spm(A')\to\spm(A)$ of sets is bijective, and via this map the affinoid subdomains of $\spm(A)$ form a basis for the strong $G$-topology on $\spm(A')$ from \cite{bgr}. 
Imposing this $G$-topology on $\spm(A)$ one gets a locally $G$-ringed space, an affinoid $k$-dagger space. (Global) $k$-dagger spaces are built from affinoid ones precisely as in \cite{bgr}. The fundamental concepts and properties from \cite{bgr} translate to $k$-dagger spaces.\\There is a faithful functor from the category of $k$-dagger spaces to the category of $k$-rigid spaces, assigning to a $k$-dagger space $X$ a $k$-rigid space $X'$ (to which we will refer as the associated rigid space). There is a natural morphism of ringed spaces $x:X'\to X$ which induces isomorphisms between the underlying $G$-topological spaces and between the stalks of the structure sheaves. $X$ is smooth if and only if $X'$ is smooth. A smooth $k$-rigid space $Y$ admits an admissible open affinoid covering $Y=\cup V_i$ such that $V_i=U_i'$ for uniquely determined (up to non canonical isomorphisms) affinoid $k$-dagger spaces $U_i$. A separated rigid space $X$ is called  partially proper, if there are  admissible open affinoid coverings $X=\cup_{j\in J} X_{j}=\cup_{j\in J} \tilde{X}_{j}$ with $X_{j}\subset\subset\tilde{X}_{j}$ for every $j\in J$ (where $\subset\subset$ is defined as in \cite{bgr}).\\
For a $k$-dagger algebra (resp. $k$-affinoid algebra) $A$, one has a universal $k$-derivation of $A$ into finite $A$-modules, $d:A\to\Omega^1_A$. In the usual way it gives rise to de Rham complexes $\Omega^{\bullet}_X=\Omega^{\bullet}_{X/k}$ on $k$-dagger (resp. $k$-rigid) spaces $X$.\\By a dagger space not specified otherwise, we will mean a $k$-dagger space, and similarly for dagger algebras, rigid spaces etc..\\
In the sequel, all dagger spaces and rigid spaces are assumed to be quasi-separated. We denote by ${\bf D}=\{x\in k; |x|\le1\}$ (resp. ${\bf D}^0=\{x\in k; |x|<1\}$) the unitdisk with (resp. without) boundary, with its canonical structure of $k$-dagger or $k$-rigid space, depending on the context. For $\epsilon\in\Gamma^*$, the ring of global functions on the polydisk $\{x\in k^n;\mbox{ all }|x_i|\le\epsilon\}$, endowed with its canonical structure of $k$-dagger space, will be denoted by $k<\epsilon^{-1}.X_1,\ldots,\epsilon^{-1}.X_n>^{\dagger}$. The dimension $\dim(X)$ of a dagger space $X$ is the maximum of all $\dim({\cal O}_{X,x})$ for $x\in X$. We say $X$ is pure dimensional if $\dim(X)=\dim({\cal O}_{X,x})$ for all $x\in X$.\\
All rigid spaces and all dagger spaces are to be understood as spaces over $k$, unless otherwise specified; they are all assumed to be quasi-separated. For a smooth {\it dagger} space $X$ we set $$R\Gamma_{dR}(X/k)=R\Gamma_{dR}(X)=R\Gamma(X,\Omega^{\bullet}_{X/k}).$$\\

\section{The Definition}

\addtocounter{satz}{1}{\bf \arabic{section}.\arabic{satz}}\newcounter{def1}\newcounter{def2}\setcounter{def1}{\value{section}}\setcounter{def2}{\value{satz}} Let $T$ be a rigid space, $S$ a dagger space with associated rigid space $S'$, and let $\phi:T\to S'$ be a closed immersion. We denote by $\Psi(\phi,S)$ the set of admissible open subsets $U$ of $S$ for which $\phi$ factors as $$T\to U'\to S'$$where $U'\to S'$ is the embedding of rigid spaces associated with the embedding of dagger spaces $U\to S$. Usually we denote the open immersion $U\to S$ by $j_U$. For an abelian sheaf ${\cal F}$ on $S$, we define $${\cal F}^{\phi}=\lim_{\stackrel{\to}{U\in\Psi(\phi,S)}}j_{U,*}{\cal F}|_U.$$

\begin{lem}\label{exa} In \arabic{def1}.\arabic{def2}, the natural morphism ${\cal F}\to{\cal F}^{\phi}$ is an epimorphism, and the functor ${\cal F}\mapsto{\cal F}^{\phi}$ is exact.
\end{lem}

{\sc Proof:} Let $S$ be affinoid and $\sigma\in\Gamma(S,{\cal F}^{\phi})$. Then $\sigma$ comes from a section $\tau\in{\cal F}(U_1)$ for some $U_1\in\Psi(\phi,S)$. Choose another $U_2\in\Psi(\phi,S)$ such that $S=(S-U_2)\cup U_1$ is an admissible covering. Obviously $\sigma|_{(S-U_2)}=0\in\Gamma(S-U_2,{\cal F}^{\phi})$ and the first claim follows. It implies the right exactness of $(?)^{\phi}$; its left exactness is clear.\\

\begin{pro}\label{basob} Let $Z$ be a rigid space, let $Y_1, Y_2$ be smooth dagger spaces, let $\phi_1:Z\to Y_1'$ be a closed immersion of $Z$ into the rigid space associated with $Y_1$, and let $\psi:Y_1\to Y_2$ be a smooth morphism of dagger spaces such that the induced morphism of rigid spaces $\psi'\circ\phi_1=\phi_2:Z\to Y_2'$ is also a closed immersion. Then the natural map$$R\Gamma(Y_2,(\Omega^{\bullet}_{Y_2})^{\phi_2})\to R\Gamma(Y_1,(\Omega^{\bullet}_{Y_1})^{\phi_1})$$is an isomorphism.
\end{pro}

{\sc Proof:} Set $Y'_Z=Y'_1\times_{Y'_2}Z$ and let $J\subset{\cal O}_{Y'_Z}$ be the ideal defining $Z$ in $Y'_Z$ (embedded diagonally). Then $J/J^2$ is locally free over ${\cal O}_Z$ since $J$ defines a section of the smooth morphism $Y'_Z\to Z$ (use \cite{EGA} 0,19.5.4). Passing to admissible coverings (the claim is local) we may suppose $Y'_Z$ is affinoid and $J/J^2$ is free. Choose a basis $\bar{t}_1,\ldots,\bar{t}_m$ of $J/J^2$, lift it to sections $t_1,\ldots, t_m\in{\cal O}_{Y'_1}$ and let $L$ be the quotient of the relative differential module $\Omega^1_{Y'_1/Y'_2}$ divided by the submodule generated by $dt_1,\ldots,dt_m$. This is a coherent ${\cal O}_{Y'_1}$-module and therefore has a Zariski closed support (\cite{bgr} 9.5.2/4). It follows from Nakayama's lemma that $L_x$ vanishes for all $x\in Y'_Z$. Thus, passing to an admissible open subset of $Y'_1$, we may suppose that $\Omega^1_{Y'_1/Y'_2}$ has $dt_1,\ldots,dt_m$ as a basis. Let $W$ be the Zariski closed subspace of $Y'_1$ defined by the ideal $(t_1,\ldots,t_m)\subset{\cal O}_{Y'_1}$. As in \cite{EGA} IV 17.4.1, 17.6.1, 17.12.1, we see that $W\to Y'_2$ is \'{e}tale. Covering admissibly and shrinking we may, in view of \ref{marklem} below, even suppose that $W\to Y_2'$ is an isomorphism, that is, we may suppose $\psi'$ has a section $\sigma:Y_2'\to Y_1'$, compatible with $\phi_1$ and $\phi_2$. We may suppose $Y_2$ and $Y_1$ are affinoid, and by \cite{kidr}, 1.18 we may suppose that there is an isomorphism $$\alpha:Y'_1\cong Y'_2\times\spm(k<\delta^{-1}.T_1,\ldots,\delta^{-1}.T_m>)$$ for some $\delta\in \Gamma^*$, where the section $\sigma$ on the left hand side corresponds to the zero section on the right hand side. For affinoid $U\in\Psi(\phi_2,Y_2)$ and $0<\epsilon\le\delta$ let $Y_1(\epsilon,U)\subset Y_1$ be the open dagger subspace of $Y_1$ defined by the open rigid subspace $$Y'_1(\epsilon,U)=\alpha^{-1}(U'\times\spm(k<\epsilon^{-1}.T_1,\ldots,\epsilon^{-1}.T_m>))$$ of $Y'_1$. Given such a $U$ we have $\psi(Y_1(\epsilon,U))\subset U$ for sufficiently small $\epsilon$. If $U$ is a Weierstrass domain in $Y_2$, then $Y_1'(\epsilon,U)$ is a Weierstrass domain in $Y_1'$, so $Y_1(\epsilon,U)$ is a Weierstrass domain in $Y_1$ (if necessary, modify the defining functions slightly to get overconvergent ones). In particular, $Y_1(\epsilon,U)$ is affinoid, hence $$Rj_{Y_1(\epsilon,U),*}\Omega^{\bullet}_{Y_1(\epsilon,U)}=j_{Y_1(\epsilon,U),*}\Omega^{\bullet}_{Y_1(\epsilon,U)}.$$     
The set of all such $Y_1(\epsilon,U)$ with $0<\epsilon\le\delta$ and Weierstrass domains $U\in\Psi(\phi_2,Y_2)$ is cofinal in the system $\Psi(\phi_1,Y_1)$. Since $Y_1$ and $Y_2$ are quasi-compact, cohomology commutes with the direct limit, so it is now enough to show that $\psi|_{Y_1(\epsilon,U)}$ induces an isomorphism $$R\Gamma_{dR}(U)\to R\Gamma_{dR}(Y_1(\epsilon,U))$$ for all such $U, \epsilon$. By \cite{boartin} we can find a map of dagger spaces $\eta:U\to Y_1(\epsilon/2,U)$ whose completion (on the level of algebras) is close to the map of rigid spaces $\sigma|_{U'}:U'\to Y'_1(\epsilon/2,U)$, namely so close that $(\psi|_{Y_1(\epsilon,U)})\circ\eta$ is close to the identity. In particular $(\psi|_{Y_1(\epsilon,U)})\circ\eta$ is an automorphism of $U$, thus induces an automorphism of $R\Gamma_{dR}(U)$. Therefore it is enough to see that $\eta$ induces an isomorphism in de Rham cohomology. Let $\delta:Y_1'(\epsilon,U)\to U'$ be the composition (of the respective restrictions) of $\alpha$ with the projection onto the first factor in the target of $\alpha$. By \cite{boartin} we can approximate $\delta$ by a map of dagger spaces $$\gamma:Y_1(\epsilon,U)\to U$$ such that $\gamma'\circ\sigma:U'\to U'$ is close to the identity, hence $\gamma\circ\eta:U\to U$ is close to the identity. In particular $\gamma\circ\eta$ is an automorphism and induces an automorphism of $R\Gamma_{dR}(U)$. So we only need to show that $\gamma$ induces an isomorphism in de Rham cohomology. We can find $S_i\in{\cal O}_{Y_1(\epsilon,U)}(Y_1(\epsilon,U))$ close to $\alpha^*(T_i)\in{\cal O}_{Y'_1(\epsilon,U)}(Y'_1(\epsilon,U))$ such that the rule $\epsilon^{-1}.T_i\mapsto\epsilon^{-1}.S_i$ defines an extension $$\tilde{\gamma}:Y_1(\epsilon,U)\to U\times\spm(k<\epsilon^{-1}.T_1,\ldots,\epsilon^{-1}.T_m>^{\dagger})$$
of $\gamma$ and such that the completion $\tilde{\gamma}'$ of $\tilde{\gamma}$ is close to the map induced by $\alpha$. In particular $\tilde{\gamma}'$ and hence $\tilde{\gamma}$ is an isomorphism. But then $\gamma$ must induce an isomorphism in de Rham cohomology.\\ 

\begin{lem}\label{marklem} Let $Z\stackrel{i}{\to}X\stackrel{f}{\to}Y$ be morphisms of affinoid rigid spaces such that $f$ is \'{e}tale and $i$ and $j=f\circ i$ are closed immersions. Then there is an admissible covering $Y=\cup_{i\in I}Y_i$ and for all $i\in I$ there are open neighbourhoods $U_i\subset Y_i$ of $Y_i\cap j(Z)$ and $V_i\subset X_i=f^{-1}(Y_i)$ of $X_i\cap i(Z)$ such that $f$ induces isomorphisms $V_i\cong U_i$.
\end{lem}

{\sc Proof:} According to \cite{djvdp} 3.1.4 we may assume, after passing to an admissible covering of $Y$, that there is a finite \'{e}tale morphism $\bar{f}:\bar{X}\to Y$ and an open immersion $l:X\to\bar{X}$ for some rigid space $\bar{X}$ such that $f=\bar{f}\circ l$. By \cite{kisin} 2.5 we may assume $Z$ is connected. Since $\bar{f}$ is \'{e}tale, there is a decomposition $\bar{f}^{-1}(j(Z))=i(Z)\coprod Q$. Again by \cite{kisin} 2.5 we find an admissible open $T\subset\bar{X}$ which decomposes as $T=T_1\coprod T_2$ with $i(Z)\subset T_1$ and $Q\subset T_2$. 
By  \cite{kisin} 2.4 there is an open connected $S\subset Y$ with $j(Z)\subset S$ and $\bar{f}^{-1}(S)\subset T$. Let $W$ be the connected component of $\bar{f}^{-1}(S)\cap T_1$ which contains $i(Z)$. Then ${\cal O}_W$ is a locally free finite rank ${\cal O}_S$-module --- its rank is one since this is so modulo the ideal defining $Z$ in $S$. Therefore $W\to S$ is an isomorphism, and $V=W\cap X\subset X$ and $U=f(V)\subset Y$ do the job.\\

\addtocounter{satz}{1}{\bf \arabic{section}.\arabic{satz}}\newcounter{df1}\newcounter{df2}\setcounter{df1}{\value{section}}\setcounter{df2}{\value{satz}} Let $T,S,\phi$ be as in \arabic{def1}.\arabic{def2} and suppose $S$ is smooth. Let $s:S'\to S$ be the natural morphism of ringed spaces. Define the de Rham cohomology of $T$ by$$R\Gamma_{dR}(T/k)=R\Gamma_{dR}(T)=R\Gamma(T,(s\circ\phi)^{-1}\Omega^{\bullet}_{S/k})$$and let $H_{dR}^q(T)=H^q(R\Gamma_{dR}(T))$. Clearly $$R\Gamma_{dR}(T)=R\Gamma(S,(\Omega^{\bullet}_{S/k})^{\phi})$$and it is this formulation to which we refer in our proof of the welldefinedness of $R\Gamma_{dR}(T/k)$.\\

\begin{pro}\label{welldef} $R\Gamma_{dR}(T)$ is independent of $S$ and $\psi$; it depends only on the reduced structure of $T$. The de Rham cohomology is a contravariant functor in $T$, with the following property: Given  $T_1,S_1,\phi_1$ and  $T_2,S_2,\phi_2$ as above and morphisms $\beta:S_1\to S_2$ and $\gamma:T_1\to T_2$ such that for the map $\beta':S'_1\to S'_2$ associated with $\beta$ we have $\beta'\circ\phi_1=\phi_2\circ\gamma$, then the map of functoriality $R\Gamma_{dR}(T_2)\to R\Gamma_{dR}(T_1)$ induced by $\gamma$ is the natural map $$R\Gamma(S_2,(\Omega^{\bullet}_{S_2})^{\phi_2})\to R\Gamma(S_1,(\Omega^{\bullet}_{S_1})^{\phi_1})$$induced by $\beta$.
\end{pro}

{\sc Proof:} Let $S_1, S_2$ be smooth dagger spaces with associated rigid spaces $S_1', S_2'$, and let $S_1'\stackrel{\phi_1}{\leftarrow}T\stackrel{\phi_2}{\to}S'_2$ be closed immersions. We compare with the diagonal embedding $\phi_{1,2}=(\phi_1,\phi_2):T\to S'_1\times S'_2$: By \ref{basob} the projections $S_{12}=S_1\times S_2\to S_i$, which are smooth, induce isomorphisms
$$R\Gamma(S_i,(\Omega^{\bullet}_{S_i})^{\phi_i})\cong R\Gamma(S_{12},(\Omega^{\bullet}_{S_{12}})^{\phi_{12}})$$for $i=1,2$. Composing we get the wanted isomorphism$$R\Gamma(S_2,(\Omega^{\bullet}_{S_2})^{\phi_2})\cong R\Gamma(S_1,(\Omega^{\bullet}_{S_1})^{\phi_1}),$$compatible with those for a third choice $S_3, \phi_3$. Now let $T_1,S_1,\phi_1$ and $T_2,S_2,\phi_2$ be as above, and let $\gamma:T_1\to T_2$ be a morphism of rigid spaces. We have 
$$R\Gamma_{dR}(T_1)=R\Gamma(S_{12},(\Omega^{\bullet}_{S_{12}})^{(\phi_1,\phi_2\circ\gamma)})$$$$R\Gamma_{dR}(T_2)=R\Gamma(S_2,(\Omega^{\bullet}_{S_2})^{\phi_2})$$and the map of functoriality $R\Gamma_{dR}(T_2)\to R\Gamma_{dR}(T_1)$ induced by $\gamma$ is by definition the one induced from the natural projection $S_{12}=S_1\times S_2\to S_2$. Again one shows that it is independent of the $S_i, \phi_i$.\\Now let in addition $\beta:S_1\to S_2$ with $\beta'\circ\phi_1=\phi_2\circ\gamma$ be given. Then $\beta$ defines a morphism $\sigma:S_1\to S_{12}$ such that the induced morphism $S'_1\to S'_{12}$ is compatibel with the $T_1$-embeddings. 
We have to show that the morphism $$R\Gamma(S_{12},(\Omega^{\bullet}_{S_{12}})^{(\phi_1,\phi_2\circ\gamma)})\to R\Gamma(S_1,(\Omega^{\bullet}_{S_1})^{\phi_1})$$induced by $\sigma$ coincides with the isomorphism which underlies the well-definedness of $R\Gamma_{dR}(T_1)$ as described above. But this follows immediately from the fact that $\sigma$ defines a section for the canonical projection $S_1\times S_{12}\to S_1$ onto the first factor.\\

\addtocounter{satz}{1}{\bf \arabic{section}.\arabic{satz}} If the rigid space $T$ admits no immersion into a rigid space which is associated with a smooth dagger space, then to define $R\Gamma_{dR}(T)$ one can use embeddings of open pieces of $T$, as in \cite{hadr}, p.28f.\\

\addtocounter{satz}{1}{\bf \arabic{section}.\arabic{satz}}\newcounter{ex1}\newcounter{ex2}\setcounter{ex1}{\value{section}}\setcounter{ex2}{\value{satz}} Examples. (a) $H_{dR}^i(\spm(T_n))=0$ if $i>0$ and $H_{dR}^0(\spm(T_n))=k$, because this is the de Rham cohomology of the dagger space $\spm(W_n)$.\\(b) If the rigid space $X$ is partially proper and smooth, then we recover the naive definition, i.e. $$H_{dR}^*(X)=H^*(X,\Omega_X^{\bullet}).$$To see this note that $X$ is the rigid space associated to a uniquely determined smooth dagger space $\tilde{X}$ (by \cite{en1dag} 2.27), and the canonical map $$H^*(\tilde{X},\Omega_{\tilde{X}}^{\bullet})\to H^*(X,\Omega_X^{\bullet})$$is an isomorphism. Indeed, by \cite{en1dag} 3.2, the maps$$H^*(\tilde{X},\Omega_{\tilde{X}}^{i})\to H^*(X,\Omega_X^{i})$$are isomorphisms for any $i$ since $X$ is partially proper. Apply this to the morphism between the respective Hodge-de Rham spectral sequences.\\(c) For some computations of the de Rham cohomology of smooth affinoid curves and of hypersurfaces, see \cite{vdp}. They show that, in these cases, the numbers $\dim_k(H_{dR}^*(?))$ are finite and are, in fact, the "correct" Betti numbers.\\

\section{Functorial properties and some exact sequences}

\begin{pro}\label{mavi} Let $T$ be a rigid space, let $T_1, T_2$ be Zariski closed subspaces of $T$ such that $T=T_1\cup T_2$. Suppose there exists a pair $(\phi:T\to S', S)$ as before. Then there is a long exact sequence$$\ldots\to H_{dR}^q(T)\to H_{dR}^q(T_1)\oplus H_{dR}^q(T_2)\to H_{dR}^q(T_1\cap T_2)\to H_{dR}^{q+1}(T)\to\ldots.$$
\end{pro}

{\sc Proof:} Let $\phi_i:T_i\to S'$ for $i=1,2$ and $\phi_{12}:T_{12}=T_1\cap T_2\to S'$ be induced by $\phi$. It is enough to show that for every $q\ge0$ the sequence$$L=[0\to(\Omega^q_S)^{\phi}\to(\Omega^q_S)^{\phi_1}\oplus(\Omega^q_S)^{\phi_2}\to(\Omega^q_S)^{\phi_{12}}\to0]$$is exact. The claim is local on $S$. If $S$ is affinoid, the set of all $U_1\cup U_2$, resp. of all $U_i$, resp. of all $U_1\cap U_2$, where the $U_i$ run through the affinoid $U_i\in\Psi(\phi_i,S)$ for $i=1,2$, is cofinal in $\Psi(\phi,S)$, resp. in $\Psi(\phi_i,S)$, resp. in $\Psi(\phi_{12},S)$. 
By the sheaf property and the commutation of direct limits with $\Gamma(V,?)$ on quasi-compact spaces $V$, we see that $\Gamma(V,L)$ is exact on the left and in the middle for all admissible open affinoid $V\subset S$, hence $L$ is exact on the left and in the middle. Since $\Omega^q_S\to(\Omega^q_S)^{\phi_{12}}$ is an epimorphism (\ref{exa}), $L$ is also exact on the right. This implies what we want.\\

\begin{pro}\label{blasex} Let $f:T_1\to T_2$ be a quasi-compact morphism of rigid spaces. Let $Z_2$ be a Zariski closed subspace of $T_2$ and let $Z_1=Z_2\times_{T_2}T_1$. Assume that $T_1-Z_1$ maps isomorphically to $T_2-Z_2$. Assume furthermore that there exist smooth dagger spaces $S_1, S_2$, closed immersions $\phi_1:T_1\to S'_1$ and $\phi_2:T_2\to S'_2$ into their associated rigid spaces, and a quasi-compact morphism $g:S_1\to S_2$ such that the induced morphism $g':S'_1\to S'_2$ satisfies $g'\circ\psi_1=\psi_2\circ f$ and maps $S'_1-g'^{-1}(Z_2)$ isomorphically to $S'_2-Z_2$. Then there is a long exact sequence$$\ldots\to H_{dR}^q(T_2)\to H_{dR}^q(T_1)\oplus H_{dR}^q(Z_2)\to H_{dR}^q(Z_1)\to H_{dR}^{q+1}(T_2)\to\ldots.$$
\end{pro}

{\sc Proof:} First suppose $T_i=S'_i$ for $i=1,2$. Note that we have a natural transformation $$(g_*(?))^{\phi_2}\to g_*((?)^{\phi_1})$$of functors of abelian sheaves on $S_1$. We claim that $$(Rg_*{\cal F})^{\phi_2}\to Rg_*({\cal F}^{\phi_1})$$for an abelian sheaf ${\cal F}$ on $S_1$ is an isomorphism. Indeed, one has $$(R^ig_*{\cal F})^{\phi_2}=\lim_{\stackrel{\to}{U\in\Psi(\phi_2,S_2)}}j_{U,*}(R^ig_*{\cal F}|_U)=\lim_{\stackrel{\to}{U\in\Psi(\phi_2,S_2)}}R^ig_*(j_{g^{-1}(U),*}{\cal F}|_{g^{-1}(U)})
$$$$\stackrel{(1)}{\cong}R^ig_*(\lim_{\stackrel{\to}{U\in\Psi(\phi_2,S_2)}}j_{g^{-1}(U),*}{\cal F}|_{g^{-1}(U)})\stackrel{(2)}{\cong}R^ig_*(\lim_{\stackrel{\to}{U\in\Psi(\phi_1,S_1)}}j_{U,*}{\cal F}|_U)=R^ig_*({\cal F}^{\phi_1}).$$Here $(1)$ holds since $g$ is quasi-compact, and $(2)$ holds since the set of all $g^{-1}(U)$ for $U\in\Psi(\phi_2,S_2)$ is cofinal within $\Psi(\phi_1,S_1)$, also because $g$ is quasi-compact (\cite{kisin} 2.4). We obtain$$(Rg_*\Omega^{\bullet}_{S_1})^{\phi_2}\cong Rg_*((\Omega_{S_1}^{\bullet})^{\phi_1}).$$
Now let $Q$ be the mapping cone of the natural map $\Omega_{S_2}^{\bullet}\to Rg_*\Omega_{S_1}^{\bullet}$. The exactness of $(?)^{\phi_2}$ and the isomorphism just seen tell us that $Q^{\phi_2}$ is the mapping cone of $(\Omega_{S_2}^{\bullet})^{\phi_2}\to Rg_*((\Omega_{S_1}^{\bullet})^{\phi_1})$. On the other hand, our assumptions imply that $Q\to Q^{\phi_2}$ is an isomorphism, i.e. $\Omega_{S_2}^{\bullet}\to Rg_*\Omega_{S_1}^{\bullet}$ and $(\Omega_{S_2}^{\bullet})^{\psi_2}\to Rg_*((\Omega_{S_1}^{\bullet})^{\psi_1})$ have an isomorphic mapping cone. By a diagram chase according to the pattern of \cite{hadr} p.44, we conclude in this case. The general case follows formally from this and \ref{mavi}, precisely as in \cite{hadr} 4.4.\\

\begin{satz} Let $X$ be a $k$-scheme of finite type, let $H^*_{dR}(X)$ be its algebraic de Rham cohomology as defined in \cite{hadr} and let $X^{an}$ be the rigid analytification of $X$. Then there is a canonical isomorphism$$H_{dR}^*(X)=H_{dR}^*(X^{an}).$$
\end{satz}

{\sc Proof:} We assume for simplicity that there is a closed embedding $X\to Y$ into a smooth $k$-scheme $Y$. By definition, $$H^*_{dR}(X)=H^*(Y,\widehat{\Omega}_Y^{\bullet})$$ where $\widehat{\Omega}_Y^{\bullet}$ denotes the formal completion of the de Rham complex ${\Omega}^{\bullet}_Y$ on $Y$ along $X$. Similarly we may define an auxiliary de Rham cohomology theory $\widehat{H}^*_{dR}$ for rigid spaces as follows: Given a rigid space $W$, choose a closed embedding $W\to Z$ into a smooth rigid space $Z$ and let $\widehat{H}^*_{dR}(W)=H^*(Z,\widehat{\Omega}_Z^{\bullet})$ where $\widehat{\Omega}_Z^{\bullet}$ denotes the formal completion of the de Rham complex ${\Omega}^{\bullet}_Z$ on $Z$ along $W$. Based on Kiehl's extension \cite{kiend} of the theorem of formal functions to the rigid analytic context one shows just as in \cite{hadr} that this definition is independent on the choice of embedding $W\to Z$. Moreover for the resulting de Rham cohomology theory $\widehat{H}^*_{dR}$ we have just as in \cite{hadr} Propositions 4.1 and 4.4 long exact sequences for blowing up and for decomposition into Zariski closed subspaces, i.e. the analogs of Propositions \ref{mavi} and \ref{blasex} above. Now from the defintions we get natural maps$$H_{dR}^*(X)\to\widehat{H}_{dR}^*(X^{an}),$$$$H_{dR}^*(X^{an})\to\widehat{H}_{dR}^*(X^{an}).$$We claim that these maps are isomorphisms. The proof, which is the same for both maps in question, is by induction
on the dimension of $X$. First we reduce to the case where $X$ is irreducible, using Proposition \ref{mavi} resp. \cite{hadr} Proposition 4.1 and its analog for $\widehat{H}^*_{dR}$. Then we perform a resolution of singularities \cite{bier}: using Proposition \ref{blasex} resp. \cite{hadr} Proposition 4.4 and its analog for $\widehat{H}^*_{dR}$ this reduces the claim to the case where $X$ is smooth. But then it follows for $H_{dR}^*(X)\to\widehat{H}_{dR}^*(X^{an})$ from \cite{kidr} and for $H_{dR}^*(X^{an})\to\widehat{H}_{dR}^*(X^{an})$ from \arabic{ex1}.\arabic{ex2} (b) (observing that $X^{an}$ is partially proper).

\begin{pro} Let $T$ be a rigid space, let $k\subset k_1$ be a finite field extension, let $T_1=T\times_{\spm(k)}{\spm(k_1)}$. There is a canonical isomorphism$$R\Gamma_{dR}(T/k)\otimes_kk_1\cong R\Gamma_{dR}(T_1/k_1).$$
\end{pro}

{\sc Proof:} Choose $S$ and $\phi:T\to S'$ as in \arabic{df1}.\arabic{df2}. Then $S_1=S\times_{\spm(k)}{\spm(k_1)}$ is a smooth $k_1$-dagger space, and $\phi_1=(\phi\times_{\spm(k)}{\spm(k_1)}):T_1\to S'\times_{\spm(k)}{\spm(k_1)}$ is a closed immersion into its associated $k_1$-rigid space. We may regard $S_1$ also as a $k$-rigid space, and we have a map $q:S_1\to S$ of $k$-rigid spaces. It induces the wanted map
$$R\Gamma(S,(\Omega^{\bullet}_{S/k})^{\phi})\otimes_kk_1\to R\Gamma(S_1,(\Omega^{\bullet}_{S_1/k_1})^{\phi_1}).$$Since the isomorphy claim is local, we may assume $S$ and $S_1$ quasi-compact. 
Then $R\Gamma$ commutes with the direct limits. Since $\{U\times_{\spm(k)}{\spm(k_1)};\quad U\in\Psi(\phi,S)\}$ is a fundamental system in $\Psi(\phi_1,S_1)$ (use \cite{kisin} 2.4), we therefore only need to check that $$R\Gamma_{dR}(U/k)\otimes_kk_1\to R\Gamma_{dR}(U\times_{\spm(k)}{\spm(k_1)}/k_1)$$is an isomorphism for all $U\in\Psi(\phi,S)$. One can assume that $U$ is affinoid, and then it follows immediately from the cohomological acyclicity of $U$ resp. of $U\times_{\spm(k)}{\spm(k_1)}$ for coherent ${\cal O}_U$- resp. ${\cal O}_{U\times_{\spm(k)}{\spm(k_1)}}$-modules (\cite{en1dag} 3.1).\\

\begin{pro}\label{gysin} Let $f:Z\to X$ be a closed immersion of smooth pure dimensional rigid spaces, associated with a closed immersion $Z^{\dagger}\to X^{\dagger}$ of dagger spaces. With $c=\codim(f)$, there is a long exact Gysin sequence $$\ldots H^{i-2c}_{dR}(Z)\to H_{dR}^i(X)\to H_{dR}^i(X-Z)\to H_{dR}^{i-2c+1}(Z)\to\ldots.$$
\end{pro}

{\sc Proof:} This results from the corresponding Gysin sequence for $Z^{\dagger}\to X^{\dagger}$, established in \cite{en2dag} 1.16.\\

\section{Finiteness and Formal Models}
Now we assume $k=\q(R)$ for a complete discrete valuation ring $R$ of mixed characteristic, and we denote by $\bar{k}$ its residue field. For an admissible formal $\spf(R)$-scheme ${\cal X}$ with generic fibre (as rigid space) ${\cal X}_k$ and specialization map $sp:{\cal X}_k\to {\cal X}$, and a subscheme $Z\subset{\cal X}_{\bar{k}}$ of its special fibre, we denote by $]Z[_{\cal X}=sp^{-1}(Z)$ the tube of $Z$ in ${\cal X}$, an admissible open subset of ${\cal X}_k$.\\For $\epsilon\in\Gamma^*$ we denote by ${\bf D}(\epsilon)$ (resp. ${\bf D}^0(\epsilon)$) the closed (resp. open) disk of radius $\epsilon$, as rigid spaces; in particular we let ${\bf D}^0={\bf D}^0(1)$.\\ 
A rigid space $L$ is called quasi-dagger, if $L$ admits an admissible covering $L=\cup_{i\in I}L_i$ such that each $L_i$ is the rigid space associated with a dagger space. A closed immersion $N\to L$ of rigid spaces is called quasi-dagger if $L$ admits an admissible covering $L=\cup_{i\in I}L_i$ such that each $N\times_LL_i\to L_i$ is associated with a closed immersion of dagger spaces.\\Smooth rigid spaces are quasi-dagger; open subspaces of quasi-dagger spaces are quasi-dagger; analytifications of algebraic $k$-schemes are quasi-dagger.\\

\begin{satz}\label{endl} Let $L, M$ and $N$ be quasi-compact rigid spaces. Suppose we are given a quasi-dagger closed immersion $N\to L$ and an open immersion $M\to L$. If we let $T=L-(M\cup N)$, then $H_{dR}^q(T)$ is finite dimensional for all $q\in\mathbb{Z}$.
\end{satz}

{\sc Proof:} Induction on $\dim(L)$. First we use \ref{mavi} and the induction hypothesis to reduce to the case where $L$ is irreducible, and clearly we can also assume $L$ is reduced. Since $N\to L$ is quasi-dagger, we may assume, after passing to a finite covering, that $N\to L$ is associated with a closed immersion $N^0\to L^0$ into a reduced and irreducible affinoid dagger space $L^0$. As explained in \cite{en2dag} 0.1, the results of \cite{bier} 1.10, \cite{schou} imply resolution of singularities for affinoid dagger spaces. Performing a resolution of singularities in our situation, we may in view of \ref{blasex} and the induction hypothesis assume that $L^0$ is smooth. But in this case we can apply \cite{en2dag} 3.5, 3.6 to conclude.\\

\begin{kor} \label{diffabg} Let $X$ be a smooth rigid Stein space, or a smooth affinoid dagger space. All differentials $d_X^i:\Omega_X^i(X)\to\Omega_X^{i+1}(X)$ have a closed image.
\end{kor}

{\sc Proof:} See \cite{kiaub} for the definition of a Stein space. A Stein space $X$ admits (in particular) an admissible covering $X=\cup_{j\in\mathbb{N}}U_j$ by affinoid rigid spaces $U_j$ with $U_j\subset U_{j+1}$. Each $\Omega^i_{U_j}(U_j)$ is a Banach space, and we endow $\Omega_X^i(X)=\lim_{\leftarrow}\Omega_{U_j}^i(U_j)$ with the inverse limit topology, hence get a Fr\'{e}chet space. The differentials $d_X^i$ are continuous. Note that since $X$ is partially proper, and acyclic for coherent ${\cal O}_X$-modules by \cite{kiaub}, we have $$H_{dR}^i(X)=\ke(d_X^i)/\bi(d_X^{i-1})$$ by \arabic{ex1}.\arabic{ex2}. 
Now first consider the special case of smooth Stein spaces of the type$$X=\cup_{\rho'<\rho}\spm(T_n({\rho'})/I.T_n({\rho'}))$$for some fixed $\rho\in\Gamma^*$ and $I<T_n(\rho)$ such that $\spm(T_n({\rho})/I)$ is smooth. In this case, $H_{dR}^i(X)$ is finite dimensional by \ref{endl}, i.e. the image of the continuous map of Fr\'{e}chet spaces $d_X^{i-1}:\Omega_X^{i-1}(X)\to\ke(d_X^i)$ is of finite codimension, hence is closed. A general Stein space $X$ admits an admissible covering $X=\cup_{j\in\mathbb{N}}V_j$ by open Stein subspaces $V_j$ of the type just considered, and such that $V_j\subset V_{j+1}$. We claim$$H_{dR}^i(X)\stackrel{(1)}{=}\lim_{\stackrel{\leftarrow}{j}}H_{dR}^{i}(V_j)=\lim_{\stackrel{\leftarrow}{j}}(\frac{\ke(d_{V_j}^i)}{\bi(d_{V_j}^{i-1})})\stackrel{(2)}{=}\frac{\lim_{\stackrel{\leftarrow}{j}}\ke(d_{V_j}^i)}{\lim_{\stackrel{\leftarrow}{j}}\bi(d_{V_j}^{i-1})}=\frac{\ke(d_X^i)}{\lim_{\stackrel{\leftarrow}{j}}\bi(d_{V_j}^{i-1})}.$$Indeed, $(1)$ follows from \cite{ss} ch.2, Cor.5 since all $H_{dR}^{i}(V_j)$ are finite dimensional. Moreover, one has $$R^1\lim_{\stackrel{\leftarrow}{j}}\Omega^{i-1}_{V_j}(V_j)=0$$ by \cite{kiaub}, and $R^2\lim_{\stackrel{\leftarrow}{j}}\ke(d_{V_j}^{i-1})=0$ because of $\kodi(\mathbb{N})=1$. Together this implies
$R^1\lim_{\stackrel{\leftarrow}{j}}\bi(d_{V_j}^{i-1})=0$ and thus euality in $(2)$. In other words, we have $\bi(d_X^{i-1})=\lim_{\stackrel{\leftarrow}{j}}\bi(d_{V_j}^{i-1})$. Since $\Omega_X^{i}(X)=\lim_{\stackrel{\leftarrow}{j}}\Omega_{V_j}^{i}(V_j)$ is a topological isomorphism, we see that $\bi(d_X^{i-1})$ is closed in $\Omega_X^{i}(X)$ since each $\bi(d_{V_j}^{i-1})$ is closed in $\Omega_{V_j}^{i}(V_j)$.\\Now let $X$ be a smooth affinoid dagger space. Then we can find a $\delta\in\Gamma^*$ and an ideal $I\subset T_n(\delta)$ such that, if for $\delta'$ with $1<\delta'<\delta$ we set $X_{\delta'}=\spm(T_n(\delta')/I.T_n(\delta'))$, each $X_{\delta'}$ is a smooth affinoid rigid space and such that $$\Omega_X^{i}(X)=\lim_{\stackrel{\to}{1<\delta'}}\Omega^{i}_{X_{\delta'}}(X_{\delta'}).$$Via this isomorphism, we define the topology on $\Omega_X^{i}(X)$ as the direct limit topology in the category of locally $k$-convex topological vector spaces (compare \cite{en1dag} 4.2; it must not be confused with the norm topology, which is coarser). Of course, if we define the Stein spaces $X^0_{\delta'}=\cup_{\delta''<\delta}X_{\delta''}$, then it is also the direct limit topology for the isomorphism$$\Omega_X^{i}(X)=\lim_{\stackrel{\to}{1<\delta'}}\Omega^{i}_{X_{\delta'}^0}(X_{\delta'}^0),$$and in view of $\bi(d_X^{i-1})=\lim_{\stackrel{\to}{j}}\bi(d_{X_j^0}^{i-1})$ and $\ke(d_X^{i})=\lim_{\stackrel{\to}{j}}\ke(d_{X_j^0}^{i})$, our claim follows from that for Stein spaces.\\

It follows that the de Rham cohomology groups of smooth $k$-rigid Stein spaces are topologically separated for their canonical topology, hence are {\it Fr\'{e}chet spaces}. In the particular case of Drinfeld's $p$-adic symmetric spaces this has been proved directly by Schneider and Teitelbaum.\\

\begin{pro}\label{kunneth} Let $T_1, T_2$ be smooth rigid spaces. There exists a canonical isomorphism$$R\Gamma_{dR}(T_1)\otimes_k R\Gamma_{dR}(T_2)\stackrel{\cong}{\longrightarrow}R\Gamma_{dR}(T_1\times T_2).$$
\end{pro}

{\sc Proof:} One easily constructs such a map and sees that the isomorphy claim is local. Therefore we can assume $T_1$ and $T_2$ are associated with smooth dagger spaces. We then conclude by \cite{en1dag} 4.12 (which is deduced from \cite{en1dag} 4.7, i.e. our \ref{diffabg}).\\

\addtocounter{satz}{1}{\bf \arabic{section}.\arabic{satz}} Due to \ref{diffabg}, one can also derive some duality formulas from \cite{en1dag} section 4.\\

\begin{lem}\label{tubhomot} Let $f:{\cal X}\to{\cal Y}$ be an immersion of smooth formal $\spf(R)$-schemes, let $Z\to{\cal X}_{\bar{k}}$ be a closed immersion into its special fibre. The canonical map $$R\Gamma_{dR}(]Z[_{\cal Y})\to R\Gamma_{dR}(]Z[_{\cal X})$$is an isomorphism.
\end{lem}

{\sc Proof:} The claim is local, so by \cite{baci} 4.3 we may suppose $$]Z[_{\cal Y}\cong ]Z[_{\cal X}\times({\bf D}^0)^r$$such that the zero section on the right hand side corresponds to the embedding $\phi:]Z[_{\cal X}\to]Z[_{\cal Y}$ induced by $f$. We may also suppose that ${\cal Y}_k$ has a smooth underlying dagger space ${\cal Y}_k^0$; let $]Z[^{\dagger}_{\cal Y}$ be its open subspace defined by $]Z[_{\cal Y}\subset{\cal Y}_k$. Then our task is to show that$$R\Gamma_{dR}(]Z[^{\dagger}_{\cal Y})\to R\Gamma(]Z[^{\dagger}_{\cal Y},(\Omega^{\bullet}_{]Z[^{\dagger}_{\cal Y}/k})^{\phi})$$is an isomorphism. 
For an affinoid open subspace $V\subset]Z[_{\cal X}$ and $\epsilon\in\Gamma^*\cap]0,1[$ define $\tilde{V}\subset]Z[^{\dagger}_{\cal Y}$ resp. $V_{\epsilon}\subset]Z[^{\dagger}_{\cal Y}$ as the open dagger subspace corresponding to the rigid subspace $V\times({\bf D}^0)^r\subset]Z[_{\cal Y}$ resp. to $V\times({\bf D}(\epsilon))^r\subset]Z[_{\cal Y}$. It is enough to show that $$R\Gamma_{dR}(\tilde{V})\to R\Gamma(\tilde{V},(\Omega^{\bullet}_{]Z[^{\dagger}_{\cal Y}/k})^{\phi})$$
is an isomorphism for all such $V$. Now $\tilde{V}$ is quasi-compact, so $$R\Gamma(\tilde{V},(\Omega^{\bullet}_{]Z[^{\dagger}_{\cal Y}/k})^{\phi})=\lim_{\stackrel{\to}{U\in\Psi(\phi,]Z[^{\dagger}_{\cal Y})}}R\Gamma(\tilde{V},j_{U,*}\Omega^{\bullet}_{U/k})=\lim_{\stackrel{\to}{U\in\Psi(\phi|_V,\tilde{V})}}R\Gamma_{dR}(U);$$
since the $V_{\epsilon}$ are cofinal in $\Psi(\phi|_V,\tilde{V})$, we only need to show that each $$R\Gamma_{dR}(\tilde{V})\to R\Gamma_{dR}(V_{\epsilon})$$is an isomorphism. Switching back to associated rigid spaces, it is enough to show that each $$R\Gamma_{dR}(\tilde{V}')\to R\Gamma_{dR}(V_{\epsilon}')$$is an isomorphism; but this is a consequence of \ref{kunneth} and the triviality of $R\Gamma_{dR}(({\bf D}^0)^r)$ and $R\Gamma_{dR}(({\bf D}(\epsilon))^r)$.\\

\begin{pro}\label{rigco} Let ${\cal X}$ be a smooth formal $\spf(R)$-scheme and let $Z\to{\cal X}_{\bar{k}}$ be a closed immersion. There is for all $q\in\mathbb{Z}$ a canonical isomorphism$$\alpha:H_{rig}^q(Z)\cong H_{dR}^q(]Z[_{\cal X}).$$If ${\cal X}_i$ for $i=1,2$ are two smooth formal $\spf(R)$-schemes with closed immersions $Z_i\to({\cal X}_i)_{\bar{k}}$, and if there is a morphism ${\cal X}_1\to{\cal X}_2$ covering a morphism $\rho:Z_1\to Z_2$, in particular inducing a morphism $\sigma:]Z_1[_{{\cal X}_1}\to]Z_2[_{{\cal X}_2}$, then $\alpha_1\circ H_{rig}^q(\rho)=H_{dR}^q(\sigma)\circ\alpha_2$.
\end{pro}

{\sc Proof:} First assume there is an embedding ${\cal X}\to{\cal Y}$ into a smooth proper $\spf(R)$-scheme ${\cal Y}$. By \cite{en1dag} 2.27 the proper rigid space ${\cal Y}_k$ has a unique underlying dagger space ${\cal Y}^{\dagger}_k$; denote by $]Z[^{\dagger}_{\cal Y}$ its admissible open subset corresponding to $]Z[_{\cal Y}\subset {\cal Y}_k$. Then $H_{dR}^q(]Z[_{\cal X})=H_{dR}^q(]Z[_{\cal Y})$ by \ref{tubhomot}, and $H_{dR}^q(]Z[_{\cal Y})=H_{dR}^q(]Z[^{\dagger}_{\cal Y})$ by definition of $H_{dR}^q(]Z[_{\cal Y})$. 
Thus our wanted isomorphism is obtained from the isomorphism$$H_{rig}^q(Z)\cong H_{dR}^q(]Z[^{\dagger}_{\cal Y})$$we constructed in \cite{en1dag} 5.1. In general, if there is no global embeding ${\cal X}\to{\cal Y}$ as above, one covers ${\cal X}$ by open affine formal subschemes and works with simplicial formal schemes and analytic spaces as usual. The functoriality assertion follows similarly as in \ref{welldef}.\\

From \ref{rigco} it follows in particular that the rigid cohomology groups of algebraic $\bar{k}$-schemes depend functorially on their rigid analytic tubes in arbitrary smooth formal $\spf(R)$-schemes, i.e. morphisms between these tubes are enough to define morphism in rigid cohomology; no extension requirements are needed.\\

\begin{pro}\label{homot} Let $X=\spm(A)$ be a smooth affinoid rigid space, let $f_1,\ldots,f_m$ be elements of $A$, and let $\epsilon\in \Gamma^*$. Set $Z=V(f_1,\ldots,f_m)$ and $U_i=\spm(A<(\epsilon^{-1}.f_i)^{-1}>$, and then $U=\cup_{i=1}^mU_i$. Suppose there is an isomorphism$$X-U\cong({\bf D}^0(\epsilon))^m\times Z$$where the functions $f_1,\ldots,f_m$ on the left hand side correspond to the standard coordinates on the right hand side. Then the natural homomorphism$$R\Gamma_{dR}(X-Z)\to R\Gamma_{dR}(U)$$is an ismorphism.
\end{pro}

{\sc Proof:} For $\delta\in\Gamma^*, \delta\le\epsilon$, and $1\le i\le m$ let $$V_{\delta}=\spm(A<\delta^{-1}.f_1,\ldots,\delta^{-1}.f_m>)$$$$U_{{\delta},i}=\spm(A<(\delta^{-1}.f_i)^{-1}>)$$and let $U_{\delta}=\cup_{i=1}^m U_{{\delta},i}$. Since for each $\delta$, the covering $X=U_{\delta}\cup V_{\delta}$ is admissible, one easily sees that it suffices to prove that the natural maps \begin{gather} \lim_{\stackrel{\to}{\delta<\epsilon}}R\Gamma_{dR}(U_{\delta})\to R\Gamma_{dR}(U)\tag{$1$}\\ R\Gamma_{dR}(V_{\delta}-Z)\to R\Gamma_{dR}(V_{\delta}\cap U_{\delta})\quad\quad (\delta<\epsilon) \tag{$2$}\end{gather} are isomorphisms. 
Since $X$ is smooth, it is quasi-algebraic (\cite{en1dag} 2.18), hence we may assume it is associated with a smooth affinoid dagger space $X^{\dagger}$. For $\delta\le\epsilon$ let $U_{{\delta},i}^{\dagger}$, resp. $U_{\delta}^{\dagger}$ be the open subspace of $X^{\dagger}$ corresponding to $U_{{\delta},i}$, resp. $U_{\delta}$. 
Then $$R\Gamma_{dR}(U)=R\Gamma_{dR}(U_{\epsilon}^{\dagger})\quad\mbox{ and }\quad R\Gamma_{dR}(U_{\delta})=R\Gamma_{dR}(U_{\delta}^{\dagger}),$$and we have admissible coverings $U_{\delta}^{\dagger}=\cup_{i=1}^mU_{{\delta},i}^{\dagger}$ for all $\delta\le\epsilon$. By the exactness of direct limits, that $(1)$ is an isomorphism will follow once we know that for all $I\subset\{1,\ldots,m\}$, the map $$\lim_{\stackrel{\to}{\delta<\epsilon}}R\Gamma_{dR}(\cap_{i\in I}U^{\dagger}_{{\delta},i})\to R\Gamma_{dR}(\cap_{i\in I}U^{\dagger}_{{\epsilon},i})$$is an isomorphism. But since all $\cap_{i\in I}U^{\dagger}_{{\delta},i}$ for $\delta\le\epsilon$ are affinoid, this follows from the very definition of a dagger algebra. Now fix $\delta<\epsilon$. Observe that $Z$ is a smooth rigid space since its fibre product with $({\bf D}^0(\epsilon))^m$ is smooth, therefore we find as above an affinoid dagger space $Z^{\dagger}$ such that $Z$ is the associated rigid space of $Z^{\dagger}$. 
Thus we get an isomorphism of the rigid space $({\bf D}^0(\epsilon))^m\times Z$ and hence of the rigid space $X-U$ with the associated rigid space of the dagger space $({\bf D}^0(\epsilon)^{\dagger})^m\times Z^{\dagger}$ (by ${\bf D}^0(\epsilon)^{\dagger}$ we mean the open disk of radius $\epsilon$ with its structure of dagger space). Under this isomorphism, the open immersion $(V_{\delta}\cap U_{\delta})\to(V_{\delta}-Z)$ is the one associated with the open immersion of dagger spaces $W_{\delta}^{\dagger}\to Y_{\delta}^{\dagger}$, where we defined
$$Y_{\delta}^{\dagger}=(\spm(k<\delta^{-1}.T_1,\ldots,\delta^{-1}.T_m>^{\dagger})-\{0\})\times Z^{\dagger},$$$$W_{{\delta},i}^{\dagger}=(\spm(k<\delta^{-1}.T_1,\ldots,\delta^{-1}.T_m,(\delta^{-1}.T_i)^{-1}>^{\dagger}))\times Z^{\dagger}$$
and $W_{\delta}^{\dagger}=\cup_{i=1}^mW_{{\delta},i}^{\dagger}$, regarded as open subspaces of $({\bf D}^0(\epsilon)^{\dagger})^m\times Z^{\dagger}$. Therefore, to prove the isomorphy of $(2)$, we only need to prove the isomorphy of$$R\Gamma_{dR}(Y^{\dagger}_{\delta})\to R\Gamma_{dR}(W^{\dagger}_{\delta}).$$
There is a natural morphism from the Cech spectral sequence computing $R\Gamma_{dR}(Y^{\dagger}_{\delta})$ by means of the admissible covering $Y_{\delta}^{\dagger}=\cup_{i=1}^m(Y_{\delta}^{\dagger}-V(T_i))$, to the Cech spectral sequence computing $R\Gamma_{dR}(W^{\dagger}_{\delta})$ by means of the admissible covering $W_{\delta}^{\dagger}=\cup_{i=1}^mW_{{\delta},i}^{\dagger}$. It shows that we only need to prove that $$R\Gamma_{dR}(\cap_{i\in I}Y^{\dagger}_{\delta}-V(T_i))\to R\Gamma_{dR}(\cap_{i\in I}W^{\dagger}_{{\delta},i})$$is an isomorphism for all $I\subset\{1,\ldots,m\}$. 
By the K\"unneth formula (\ref{kunneth}), this is reduced to proving that$$R\Gamma_{dR}(\spm(k<\delta^{-1}.T>^{\dagger})-\{0\})\to R\Gamma_{dR}(\spm(k<\delta^{-1}.T,(\delta^{-1}.T)^{-1}>^{\dagger}))$$($T$ a single variable) is an isomorphism, which results from a simple computation of both sides.\\ 

\begin{kor}\label{sshomot} (a) Let $\pi$ be a uniformizer of $R$ and let ${\cal Z}\to{\cal X}$ be a closed immersion of formal $R$-schemes. Assume that locally on ${\cal X}$ there exist \'{e}tale morphisms of formal $R$-schemes$$q:{\cal X}\to\spf(R<X_1,\ldots,X_n>/(X_1\ldots X_r-\pi))$$such that ${\cal Z}=\cap_{j=r+1}^mV(q^*X_j)$ for some $1\le r\le m\le n$. Then the natural map$$R\Gamma_{dR}({\cal X}_k-{\cal Z}_k)\to R\Gamma_{dR}(]{\cal X}_{\bar{k}}-{\cal Z}_{\bar{k}}[_{\cal X})$$is an isomorphism.\\(b) In (a), the canonical map $R\Gamma_{dR}(]{\cal Z}_{\bar{k}}[_{\cal X})\to R\Gamma_{dR}({\cal Z}_k)$ is an isomorphism.\\(c) In the description of the setting in (a), replace ${\cal Z}=\cap_{j=r+1}^mV(q^*X_j)$ by ${\cal Z}=\cup_{j=r+1}^mV(q^*X_j)$. Then the same statement as in (a) holds.
\end{kor}

{\sc Proof:} (a) Since the assertion is local on ${\cal X}$, we may assume there exists globally a map $q$ as described. By \cite{en2dag} 2.6, there exists an isomorphism \begin{gather} ]{\cal Z}_{\bar{k}}[_{\cal X}\cong {\cal Z}_k\times({\bf D}^0)^{m-r}\tag{$*$}\end{gather} where the $q^*X_j$ for $r+1\le j\le m$ correspond to the standard coordinates. Thus for (a) we can cite \ref{homot}. Also (b) follows from $(*)$. Assertion (c) follows formally from (a) by Cech complex arguments.\\

\begin{kor}\label{rigcogysin} In \ref{sshomot} (a), assume that the closed embedding ${\cal Z}_k\to {\cal X}_k$ of rigid spaces is associated with a closed embedding of dagger spaces. Then there is a long exact Gysin sequence $$\ldots\to H^{i-2(m-r)}_{dR}(]{\cal Z}_{\bar{k}}[_{\cal X})\to H_{dR}^i({\cal X}_k)\to H_{dR}^i(]{\cal X}_{\bar{k}}-{\cal Z}_{\bar{k}}[_{\cal X})\to H_{dR}^{i-2(m-r)+1}(]{\cal Z}_{\bar{k}}[_{\cal X})\to\ldots.$$ 
\end{kor}

{\sc Proof:} Combine \ref{sshomot} with \ref{gysin}.\\

{\bf Remark:} (1) In \ref{rigcogysin}, the assumption that ${\cal Z}_k\to {\cal X}_k$ is associated with a closed embedding of dagger spaces, is fulfilled in particular if ${\cal Z}\to {\cal X}$ is obtained by $\pi$-adic completion from a closed immersion of $R$-schemes of finite type, or if ${\cal X}$ is proper.\\
(2) Suppose $r=1$ in \ref{rigcogysin}. Then ${\cal X}$ and ${\cal Z}$ are smooth, and we recover the long exact Gysin sequence$$\ldots\to H^{i-2(m-r)}_{rig}({\cal Z}_{\bar{k}})\to H_{rig}^i({\cal X}_{\bar{k}})\to H_{rig}^i({\cal X}_{\bar{k}}-{\cal Z}_{\bar{k}})\to H_{rig}^{i-2(m-r)+1}({\cal Z}_{\bar{k}})\to\ldots$$for rigid cohomology, constructed in \cite{berfi} and \cite{mebend}. For general $r$, our sequence might be thought of as its version for logarithmic rigid cohomology, where the base $\spf(R)$ is endowed with its canonical log structure.\\
(3) From \ref{gysin} and \ref{homot} one can derive more general Gysin sequences, without knowledge of reductions. Such occured for example in \cite{desha} p. 186, p. 190.\\

\section{Vanishing cycles}

Let ${\cal X}$ be an admissible proper formal $\spf(R)$-scheme. Attached to it are two de Rham type cohomology theories, finite dimensional over $k$: The rigid cohomology $R\Gamma_{rig}({\cal X}_{\bar{k}}/k)$ of its special fibre ${\cal X}_{\bar{k}}$, and the de Rham cohomology $R\Gamma_{dR}({\cal X}_k)$ of its generic fibre (as rigid space) ${\cal X}_k$. We want to compare these two.\\For simplicity we assume that there exists a smooth admissible formal $\spf(R)$-scheme ${\cal Y}$ and an embedding ${\cal X}\to{\cal Y}$. It induces an embdding ${\cal X}_k\to{\cal Y}_k$ into the generic fibre (as rigid space) ${\cal Y}_k$ of ${\cal Y}$. Since ${\cal X}_k$ is proper we do not need dagger spaces to define $R\Gamma_{dR}({\cal X}_k)$: We have$$R\Gamma_{dR}({\cal X}_k)=R\Gamma(]{\cal X}_{\bar{k}}[_{\cal Y},\lim_{\stackrel{\to}{U}}j_{U,*}\Omega_U^{\bullet})$$where $j_U:U\to ]{\cal X}_{\bar{k}}[_{\cal Y}$ runs through the inclusions of admissible open subsets of $]{\cal X}_{\bar{k}}[_{\cal Y}$ containing ${\cal X}_k$. On the other hand we have$$R\Gamma_{rig}({\cal X}_{\bar{k}}/k)=R\Gamma(]{\cal X}_{\bar{k}}[_{\cal Y},\Omega_{]{\cal X}_{\bar{k}}[_{\cal Y}}^{\bullet}).$$Denote by $s:]{\cal X}_{\bar{k}}[_{\cal Y}\to{\cal X}_{\bar{k}}$ the specialization map. Note that $$s_*\lim_{\stackrel{\to}{U}}j_{U,*}\Omega_U^{\bullet}=Rs_*\lim_{\stackrel{\to}{U}}j_{U,*}\Omega_U^{\bullet}$$$$s_*\Omega_{]{\cal X}_{\bar{k}}[_{\cal Y}}^{\bullet}=Rs_*\Omega_{]{\cal X}_{\bar{k}}[_{\cal Y}}^{\bullet}.$$This follows from the acyclicity of coherent modules on quasi-Stein spaces \cite{kiaub}. Thus$$R\Gamma_{dR}({\cal X}_k)=R\Gamma({\cal X}_{\bar{k}},s_*\lim_{\stackrel{\to}{U}}j_{U,*}\Omega_U^{\bullet})$$$$R\Gamma_{rig}({\cal X}_{\bar{k}}/k)=R\Gamma({\cal X}_{\bar{k}},s_*\Omega_{]{\cal X}_{\bar{k}}[_{\cal Y}}^{\bullet}).$$The canonical restriction map of sheaf complexes on ${\cal X}_{\bar{k}}$ $$\tau:s_*\Omega_{]{\cal X}_{\bar{k}}[_{\cal Y}}^{\bullet}\to s_*\lim_{\stackrel{\to}{U}}j_{U,*}\Omega_U^{\bullet}$$is injective. In other words, we can filter the sheaf complex $s_*\lim_{\stackrel{\to}{U}}j_{U,*}\Omega_U^{\bullet}$ on ${\cal X}_{\bar{k}}$ which computes $R\Gamma_{dR}({\cal X}_k)$ by a subcomplex $\bi(\tau)$ which computes $R\Gamma_{rig}({\cal X}_{\bar{k}}/k)$.\\The situation is similar to that in Bruno Chiarellotto's paper \cite{chia}: There, for a proper semistable $\bar{k}$-log scheme $Y$ he shows that the Hyodo-Steenbrink complex $WA^{\bullet}$ which computes the Hyodo-Kato cohomology of $Y$ --- which is isogenous to the de Rham cohomology of a semistable lift of $Y$, if such a lift exists --- can be filtered by a subcomplex which computes the rigid cohomology of $Y$, provided the $p$-adic monodromy-weight conjecture holds for $Y$.\\Returning to our situation, define the sheaf complex ${\cal V}^{\bullet}_{{\cal X},{\cal Y}}$ on ${\cal X}_{\bar{k}}$ by the exact sequence$$0\to s_*\Omega_{]{\cal X}_{\bar{k}}[_{\cal Y}}^{\bullet}\stackrel{\tau}{\to}s_*\lim_{\stackrel{\to}{U}}j_{U,*}\Omega_U^{\bullet}\to{\cal V}^{\bullet}_{{\cal X},{\cal Y}}\to0.$$One might call it the complex of vanishing cycles for ${\cal X}$ with respect to ${\cal Y}$: It measures the difference between $R\Gamma_{dR}({\cal X}_k)$ and $R\Gamma_{rig}({\cal X}_{\bar{k}}/k)$; we have $R\Gamma_{dR}({\cal X}_k)=R\Gamma_{rig}({\cal X}_{\bar{k}}/k)$ if $H^q({\cal X}_{\bar{k}},{\cal V}^{\bullet}_{{\cal X},{\cal Y}})=0$ for {\it all} $q\ge0$ (somewhat different to the $l$-adic vanishing cycles). For example, if ${\cal X}_k$ is smooth, or more specifically if ${\cal X}$ is the formal completion of a proper $R$-scheme of finite type with smooth generic fibre, this might be a starting point for a rigid analytic investigation of the monodromy on $H^*_{dR}({\cal X}_k)$ due to bad reduction (note that even if ${\cal X}_k$ is smooth there is no obvious subcomplex of the sheaf complex $\Omega^{\bullet}_{{\cal X}_k}$ on ${\cal X}_k$ which computes $R\Gamma_{rig}({\cal X}_{\bar{k}}/k)$; so even if ${\cal X}_k$ is smooth the construction of this paper has its assets).\\

{\it Acknowledgments: I wish to thank R. Huber for his interest and for suggesting to formulate the definition of $H^*_{dR}(X)$ as above (in contrast to the description which will be used in proof of welldefinedness ). I thank Mark Kisin for pointing out that the proof given in \cite{en1dag} 4.7 is incomplete; here we fill the gap, see \ref{diffabg}.} \\


\end{document}